\documentclass[11pt]{article}
\usepackage{amsfonts}
\usepackage{lscape}
\usepackage{latexsym,amsmath,amsthm}
\usepackage{amssymb,array}
\graphicspath{{figures/}}
\usepackage{subfigure}
\usepackage{fancyhdr}
\usepackage{multicol}
\usepackage{float}
\usepackage{csquotes}
\makeatletter 
\RequirePackage[dvips]{graphicx} 
\newtheorem{t1}{Theorem}[section]
\newtheorem{l1}{Lemma}[section]

\usepackage[numbers,round]{natbib}
\usepackage{url}

\begin{document}
\title{\textbf{On estimation of the PMF and the CDF of a natural discrete one parameter polynomial exponential distribution}}
\author{ Indrani Mukherjee$^1$, Sudhansu S. Maiti$^1$\footnote{Corresponding author e-mail:
dssm1@rediffmail.com} and Rama Shanker$^2$ \\
$^1$Department of Statistics, Visva-Bharati
University, Santiniketan-731 235, India\\
$^2$ Department of Statistics, Assam University, Silchar-7880111, India}
\date{}
\maketitle
\begin{center}
Abstract
\end{center}
In this article, a new natural discrete analog of the one parameter polynomial exponential (OPPE) distribution as a mixture of a number of negative binomial distributions has been proposed and is called as a natural discrete one parameter polynomial exponential (NDOPPE) distribution. This distribution is a generalized version of natural discrete Lindley (NDL) distribution, proposed and studied by Ahmed and Afify (2019). Two estimators viz., MLE and UMVUE of the PMF and the CDF of a NDOPPE distribution have been derived. The estimators have been compared with respect to their MSEs. Simulation study has been conducted to verify the consistency of the estimators. A real data illustration has been reported.

\vspace{0.5cm}


\section{Introduction}
\ Some standard discrete distributions have been mentioned and the estimators of their probability mass functions (PMF) and cumulative distribution functions (CDF) are studied in Maiti and Mukherjee (2017). Sometimes these models are too restrictive. As for example, Poisson model is not appropriate because it imposes the restriction of equidispersion in the modeled data. Similarly, binomial model imposes the restriction of under-dispersion. As a result, various models prescribed using discrete concentration and discrete analog approaches that are less restrictive (see, Nakagawa and Osaki (1975), Famoye (1993), among others). The most recent discrete distributions are due to Stein and Dattero(1984), Roy (2002), Roy (2003), Roy (2004), Krishna and Pundir (2009), Jazi et al. (2010), Gomez-Deniz (2010), Gomez-Deniz and Calderin-Ojeda (2011), Maiti et al. (2018), among others.

\par To generate a natural discrete analog of the one parameter polynomial exponential (OPPE)  distribution, we use the well-known fact that the geometric and the negative binomial distributions are the natural discrete analogs of the exponential and the gamma distributions respectively. This is the motivation behind proposing a new natural discrete analog to the OPPE distribution by mixing discrete counterparts of the exponential and gamma distributions.

\par The PMF of the random variable $X$ of a natural discrete one parameter polynomial exponential distribution (NDOPPE) distribution can be written as
\begin{equation}
f(x)=h(\theta)p(x)(1-\theta)^x, ~ x=0,1,\ldots,~ 0<\theta<1,
\end{equation}
where, $h(\theta)=\frac{1}{\sum_{k=1}^{r}a_{k-1}\frac{\Gamma(k)}{\theta^{k}}}$ ,
$p(x)=\sum_{k=1}^{r}a_{k-1}\Gamma(k)\binom{x+k-1}x.$

\vspace{0.1in}
The distribution can also be written as
\begin{equation}\label{ch6doppepmf}
f(x)=h(\theta) \sum_{k=1}^r a_{k-1}\frac{\Gamma(k)}{\theta^k}f_{NB}(x;\theta,k),
\end{equation}
where $f_{NB}(x;\theta,k)$ is the PMF of a negative binomial distribution with parameters $\theta$ and $k$ and $a_{k-1}$'s are known non-negative constants. 

\vspace{0.1in}
The CDF is given by
\begin{equation}\label{ch6doppecdf}
F(x)=h(\theta)\sum_{k=1}^r a_{k-1}\frac{\Gamma(k)}{\theta^k}I_\theta(x,k+1) ,~ x=0,1,\ldots,~ 
0<\theta<1,
\end{equation}
where $I_x(a,b)=\frac{1}{B(a,b)}\int_{0}^x t^{a-1}(1-t)^{b-1}dt$.

\vspace{0.2in}

\par We have mentioned two special cases which are given below,
\begin{enumerate}
	\item $r=1,~a_0=1$ gives the geometric distribution ,
	\item $r=2,~a_0=1,~a_1=1$ gives a natural discrete Lindley (NDL) distribution [c.f. Ahmed	and Afify (2019)]. The PMF and the CDF is given by
	\begin{equation*}
	f(x)=\frac{\theta^2}{(1+\theta)}(2+x)(1-\theta)^x,~ x=0,1,\ldots; \theta \in (0,1)
	\end{equation*}
	and
	\begin{equation*}
	F(x)=1-\frac{1+\theta+\theta x}{(1+\theta)}(1-\theta)^x,~ x=0,1,\ldots; \theta \in (0,1),
	\end{equation*}
	respectively.
\end{enumerate}
%We already have studied the estimation of the PMF and the CDF of geometric distribution in article 5. In this article, we have considered the second case, the NDL distribution. 
\par The problem of estimation of the PMF and the CDF is interesting for many reasons. For example, the PMF and the CDF can be used for estimation of differential entropy, R\'{e}nyi entropy, Kullback-Leibler divergence, Fisher's Information, Cumulative residual entropy, quantile function, Lorenze curve, Hazard rate function, Mean remaining life function etc.

 Most of the times, emphasis is given to infer the parameters involved in the distribution and the study is concentrated on measuring the efficiency of these estimators. No such effort has been made to find out estimator of the PMF and the CDF of these discrete random variables. Plugging in the MLE of parameter gives, by invariance property, the MLE of the PMF and the CDF. It is to be noted that all these estimators are biased. Mere substitution of the UMVUE of parameter(s) will not provide the UMVUE of the PMF and that of the CDF. Therefore, comparing UMVUE of the PMF and that of the CDF (the only unbiased estimators in our study) with other estimators seem to be an interesting study.

\par Similar type of studies have appeared in recent literature for some continuous distributions. See, Asrabadi (1990), Bagheri et al. (2014), Dixit and Jabbari (2010), Dixit and Jabbari (2011), Jabbari and Jabbari (2010), Maiti and Mukherjee (2017), Mukherjee and Maiti (2018) etc.

%Alizadeh et al.(2015), Asrabadi (1990), Bagheri et al. (2014, 2016a, 2016b), Dixit and Jabbari (2010), Dixit and Jabbari (2011), Jabbari and Jabbari (2010), Maiti and Mukherjee (2017), Mukherjee et al. (2016). 

\par The article is organized as follows. Section $\ref{ch6mlepdfcdf}$ deals with MLE of the PMF and the CDF of a NDOPPE distribution. Section $\ref{ch6umvuepdfcdf}$ is devoted to finding out the UMVUE of the PMF and the CDF. In section $\ref{ch6sim}$, simulation study results are reported and comparisons are made. Real-life data set is analyzed in section $\ref{ch6data}$. In section $\ref{ch6conclusion}$, concluding remarks are made based on the findings of this article.

\section{MLE of the PMF and that of the CDF}\label{ch6mlepdfcdf}
\ Let $X_1, X_2,...,X_n$ be a random sample of size $n$ drawn from the PMF in $\eqref{ch6doppepmf}$. here we try to find the MLE of $\theta$ which is denoted as $\widetilde{\theta}$. The log-likelihood of $\theta$ is given by
\begin{eqnarray*}
	l(\theta)&=&\ln L(\theta|X)\\
	&=&n \ln h(\theta)+\sum_{i=1}^{n}\ln p(X_{i})+\ln(1-\theta) \sum_{i=1}^{n} X_{i}\mbox{.}
\end{eqnarray*}
Now,
\begin{eqnarray}\label{ch6mle}
\frac{dl(\theta)}{d\theta}&=&0\nonumber\\
\mbox{i.e.} n \frac{d}{d\theta}\left(\ln h(\theta)\right)+\frac{1}{1-\theta}\sum_{i=1}^{n} X_{i}&=&0\nonumber\\
\mbox{i.e.} h(\theta)\sum_{k=1}^ra_{k-1}\frac{\Gamma(k+1)}{\theta^{k+1}}&=&\frac{\bar{X}}{1-\theta}\mbox{.}
\end{eqnarray}
Since, the above equation is not of a closed form, we have to solve $\eqref{ch6mle}$ numerically to obtain the MLE of $\theta$. %Using the invariance property of MLE, one can obtain the MLEs of the PMF and CDF by substituting $\widetilde{\theta}$ in $\eqref{ch6doppepmf}$ and $\eqref{ch6doppecdf}$ respectively. 
Theoretical expression for the MSE of the MLEs are not available. MSE will be studied through simulation.
\section{UMVUE of the PMF and that of the CDF}\label{ch6umvuepdfcdf}
\ In this section, we obtain the UMVUE of the PMF and that of the CDF of a NDOPPE distribution. Also, we obtain the MSEs of these estimators.

\begin{t1}\label{ch6fz}
	Let, $X_1, X_2,..., X_n \sim f_{NDOPPE}(x,\theta)$. Then the distribution of $T=X_1+X_2+....+X_n$ is
	\[
	f(t)=h^{n}(\theta)\sum_{y_1} \ldots \sum_{y_r} c(n,y_1,\ldots ,y_r)\binom{t+\sum_{k=1}^r k y_k-1}t (1-\theta)^t,
	\]
	$t=0,1,\ldots$ with $y_1+...+y_r=n$ and
	$c(n,y_1,\ldots ,y_r)=\frac{n!}{y_1!\ldots y_r!}$\\$\times \prod_{k=1}^{r}(a_{k-1}\Gamma(k))^{y_k}$.
\end{t1}
\begin{proof}
	The mgf of $T$ is
	\begin{eqnarray*}
		M_{T}(t)&=&h^n(\theta)\left[\sum_{k=1}^r \frac{a_{k-1}}{\theta ^k}\cdot \frac{\Gamma(k)\theta ^k}{\left\lbrace1-(1-\theta)e^t\right\rbrace ^k}\right]^n\\
		&=&h^n(\theta)\left[\sum_{y_1} \ldots \sum_{y_r} c(n,y_1,\ldots ,y_r)\right. \\
		&&\left.\times \frac{1}{\theta^{\sum_{k=1}^r k y_k}}\left\lbrace \frac{\theta}{1-(1-\theta)e^t}\right\rbrace^{\sum_{k=1}^rky_k} \right].
	\end{eqnarray*}
	Hence, the distribution of $T$ is
	\begin{eqnarray*}
		f(t)&=&h^{n}(\theta)\sum_{y_1} \ldots \sum_{y_r} c(n,y_1,\ldots ,y_r)\\
		&&\times \binom{t+\sum_{k=1}^r k y_k-1}t (1-\theta)^t,
	\end{eqnarray*}
	where, $c(n,y_1,\ldots ,y_r)=\frac{n!}{y_1!\ldots y_r!}\prod_{k=1}^{r}(a_{k-1}\Gamma(k))^{y_k}$.
\end{proof}

\begin{l1}\label{ch6lemma}
	The conditional distribution of $X_{1}$ given $X_{1}+X_{2}+....+X_{n}=T$ is
	\begin{eqnarray*}
		f_{X_1|T}(x|t)&=&\frac{p(x)}{A_{n}(t)} \sum_{q_1}...\sum_{y_r} c(n-1,q_1,\ldots ,q_r)\\
		&&\times \binom{t-x+\sum_{k=1}^r kq_k-1}{t-x},~x=0,\ldots ,t,
	\end{eqnarray*}
	where, 
	$$A_{n}(t)= \sum_{y_1}\ldots \sum_{y_r} c(n,~y_1,~\ldots,~y_r)\binom{t+\sum_{k=1}^r ky_k-1}t$$
	and
	$$c(n-1,q_1,\ldots,q_r)=\frac{(n-1)!}{q_1!\ldots q_r!}\prod_{k=1}^r(a_{k-1}\Gamma(k))^{y_k},$$ with $q_1+q_2+....+q_r=n-1.$
\end{l1}
\begin{proof}
	\begin{eqnarray*}
		f_{X_1|T}(x|t)&=& \frac{f(x)f(t-x)}{f(t)}\\
		&=&\frac{p(x)}{A_{n}(t)} \sum_{q_1}...\sum_{y_r} c(n-1,q_1,\ldots ,q_r)\\
		&&\binom{t-x+\sum_{k=1}^r kq_k-1}{t-x}.
	\end{eqnarray*}
\end{proof}

\begin{t1}
	Let, $T=t$ be given. Then
	\begin{eqnarray}\label{ch6umvuepmf}
	\widehat{f}(x)&=&\frac{p(x)}{A_{n}(t)} \sum_{q_1}...\sum_{q_r} c(n-1,q_1,\ldots ,q_r)\nonumber\\
	&&\times \binom{t-x+\sum_{k=1}^r kq_k-1}{t-x},~x=0,~\ldots, ~ t,
	\end{eqnarray}
	is UMVUE for $f(x)$ and
	\begin{eqnarray}\label{ch6umvuecdf}
	\widehat{F}(x)&=&\frac{1}{A_n(t)}\sum_{q_1}\ldots \sum_{q_r} c(n-1,q_1,\ldots ,q_r)\nonumber\\
	&&\times \sum_{w=0}^x p(w)\binom{t+\sum_{k=1}^r kq_k-1}{t-w},~x=0,~\ldots ,~t,
	\end{eqnarray}
	is UMVUE for $F(x)$,
	where, $p(w)=\sum_{k=1}^r a_{k-1}\Gamma(k)\binom{w+k-1}w$.
\end{t1}

\vspace{0.2in}
\textbf{Remark:}
When $r=2,~a_0=1,~a_1=1$, the above expressions $\eqref{ch6umvuepmf}$ and $\eqref{ch6umvuecdf}$ reduce to

\begin{eqnarray}
\widehat{f}(x)&=&\frac{(2+x)}{\sum_{k=0}^n \binom{n}k \binom{2n-k+t-1}t}\sum_{k=0}^{n-1} \binom{n-1}k \binom{2(n-1)-k+t-x-1}{t-x}\nonumber\\
&&,~x=0,~\ldots, ~ t
\end{eqnarray}
and
\begin{eqnarray}
\widehat{F}(x)&=&\sum_{w=0}^x \frac{(2+w)}{\sum_{k=0}^n \binom{n}k \binom{2n-k+t-1}t} \sum_{k=0}^{n-1} \binom{n-1}k \nonumber\\
&& \times \binom{2(n-1)-k+t-w-1}{t-w},~x=0,~\ldots, ~ t,
\end{eqnarray}
respectively.

The MSE of $\widehat{f}(x)$ is given by
\begin{eqnarray}\label{ch6msefhatx}
MSE(\widehat{f}(x))&=&E(\widehat{f}^2(x))-f^2(x)\nonumber\\
&=&\sum _{t=x}^\infty \left[\frac{p(x)}{A_{n}(t)}\sum_{q_0} \sum_{q_1}...\sum_{q_r} c(n-1,q_0,q_1,\ldots ,q_r)\right.\nonumber\\
&&\left.\binom{t-x+\sum_{k=1}^r kq_k-1}{t-x}\right]^2 f(t)-f^2(x).
\end{eqnarray}
Using Theorem $\ref{ch6fz}$, $\eqref{ch6doppepmf}$ in $\eqref{ch6msefhatx}$, we can get the value of the MSE of UMVUE of the PDF.

The MSE of $\widehat{F}(x)$ is given by
\begin{eqnarray}\label{ch6mseFhatx}
MSE(\widehat{F}(x))&=&E(\widehat{F}^2(x))-F^2(x)\nonumber\\
&=&\sum _{t=x}^\infty \left[\sum_{w=0}^x \frac{p(w)}{A_n(t)}\sum_{q_1}...\sum_{q_r} c(n-1,q_0,q_1,\ldots ,q_r)\right.\nonumber\\
&&\times \left.\binom{t-w+\sum_{k=1}^r kq_k-1}{t-w}\right]^2f(t)\nonumber\\
&&-F^2(x),
\end{eqnarray}
where, $p(w)=\sum_{k=1}^r a_{k-1}\Gamma(k)\binom{w+k-1}w$.
Similarly, using Theorem $\ref{ch6fz}$, $\eqref{ch6doppecdf}$ in $\eqref{ch6mseFhatx}$, we can get the value of the MSE of UMVUE of the CDF.

\par
\begin{figure}[H]
	\centering
	\includegraphics[height=3in,width=6in]{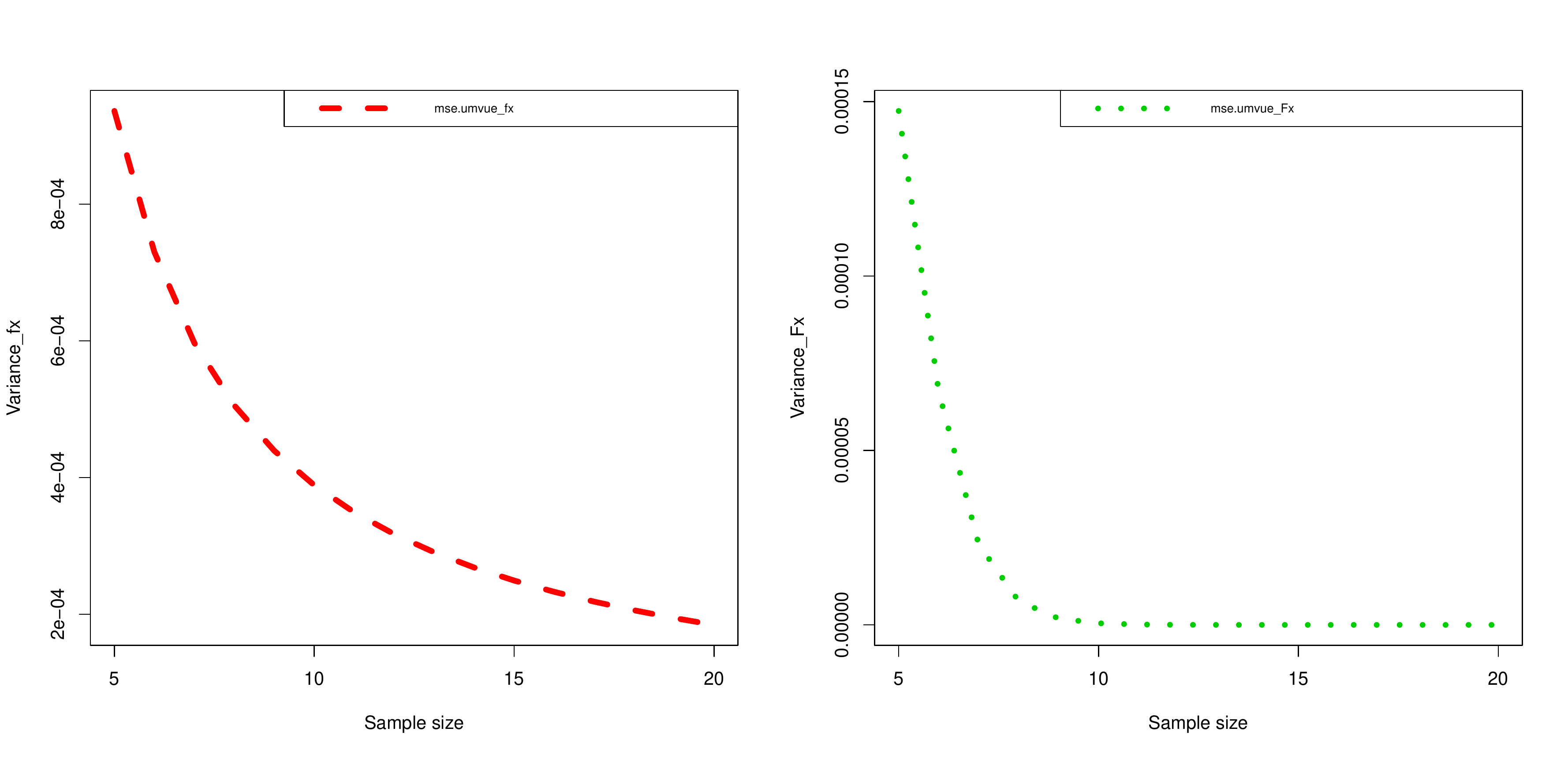}
	\caption{Graph of theoretical MSE of the UMVUE of the PMF and that of the CDF of a NDL distribution for $x=2$, $\theta=0.01$ and $r=2$.}\label{ch6theorypmfcdf}
\end{figure}

 Theoretical MSE of the UMVUE of the PMF and that of the CDF have been shown in Figure $\ref{ch6theorypmfcdf}$, taking $a_0=1,~a_1=1$, $x=2$, $\theta=0.01$ and $r=2$. It is clear from the graph that MSE (in this case, variance) decreases as sample size increases.

%We represented the theoretical graph of MSE (variance in this case) of the UMVUE of the PMF and CDF in Figure $\ref{ch6theorypmfcdf}$. From the above figure, we can say that MSE decreases with increasing sample size.
%\section{Reliability related measures}\label{ch6secrel}
%\begin{enumerate}

%\item  Reliability function: 
%\begin{equation}\label{ch6rel}
%R(x)=\sum_{j=x+1}^\infty f(j)
%\end{equation}

%\item Hazard rate function:
%\begin{equation}\label{ch6haz}
%h(x)=\frac{f(x)}{\sum_{j=x}^\infty f(j)}
%\end{equation}	

%\item Mean remaining life function: 
%\begin{equation}\label{ch6mrl}
%m(x)=\frac{1}{\sum_{j=x}^\infty f(j)} \sum_{j=x}^\infty \sum_{k=j+1}^\infty f(k)
%\end{equation}

%\end{enumerate}

%We have the value of MLE of the parameter of distribution, taking $a_0=1,~a_1=1$, $x=2$, $\theta=0.01$ and $r=2$, from Section $\ref{ch6mlepdfcdf}$, and UMVUE of the PMF, $\hat{f}(x)$ and CDF, $\hat{F}(x)$ from Section $\ref{ch6umvuepdfcdf}$ respectively. Hence plugging in the value of estimators in $\eqref{ch6rel},~\eqref{ch6haz},~\eqref{ch6mrl}$, we can easily find the corresponding estimates of reliability, hazard rate and mean remaining life function, respectively. 

\section{Simulation}\label{ch6sim}

\ Generation of random sample $X_1,~X_2,~\ldots,~X_n$ is distributed in the following algorithm:
\begin{enumerate}
	\item[1.] Generate $U_{i}\sim Uniform (0,1), i=1(1)n$.
	\item[2.] If $\frac{\sum_{k=1}^{j-1}a_{k-1}\frac{(k-1)!}{\theta ^k}}{\sum_{k=1}^r a_{k-1}\frac{(k-1)!}{\theta ^k}}< U_{i} \leq \frac{\sum_{k=1}^j a_{k-1}\frac{(k-1)!}{\theta ^k}}{\sum_{k=1}^r a_{k-1}\frac{(k-1)!}{\theta^k}}~,j=2,...,r$, then set $X_{i}=V_{i}$,\\ where $V_{i}\sim NB(j,\theta)$ and if $U_{i}\leq \frac{a_{0}\frac{1}{\theta}}{\sum_{k=1}^r a_{k-1}\frac{(k-1)!}{\theta^k}}$, then set $X_{i}=V_{i}$, where $V_{i}\sim geo(\theta).$
\end{enumerate}
%We assume that,$a_{0}=a_{1}=1,~r=2$ which leads to the discrete version of Lindley distribution. 
%In this article, we have studied the properties of the PMF and the CDF of a discrete version of Lindley distribution. In the previous article (article 5) we already study the properties of the PMF and the CDF of geometric distribution. 
\par
\begin{figure}[H]
	\centering
	\subfigure{\includegraphics[height=3in,width=2.6in]{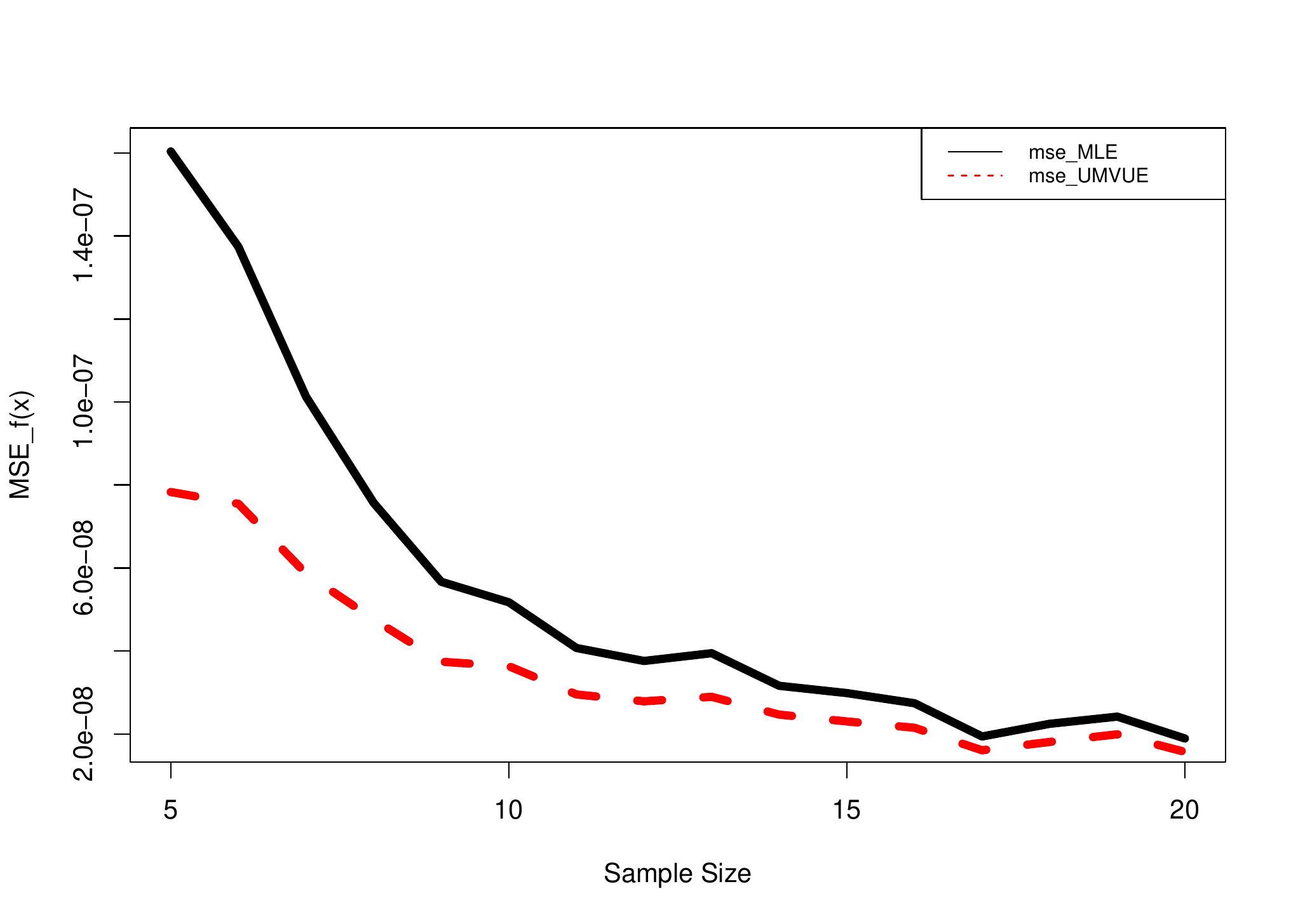}}
	\subfigure{\includegraphics[height=3in,width=2.6in]{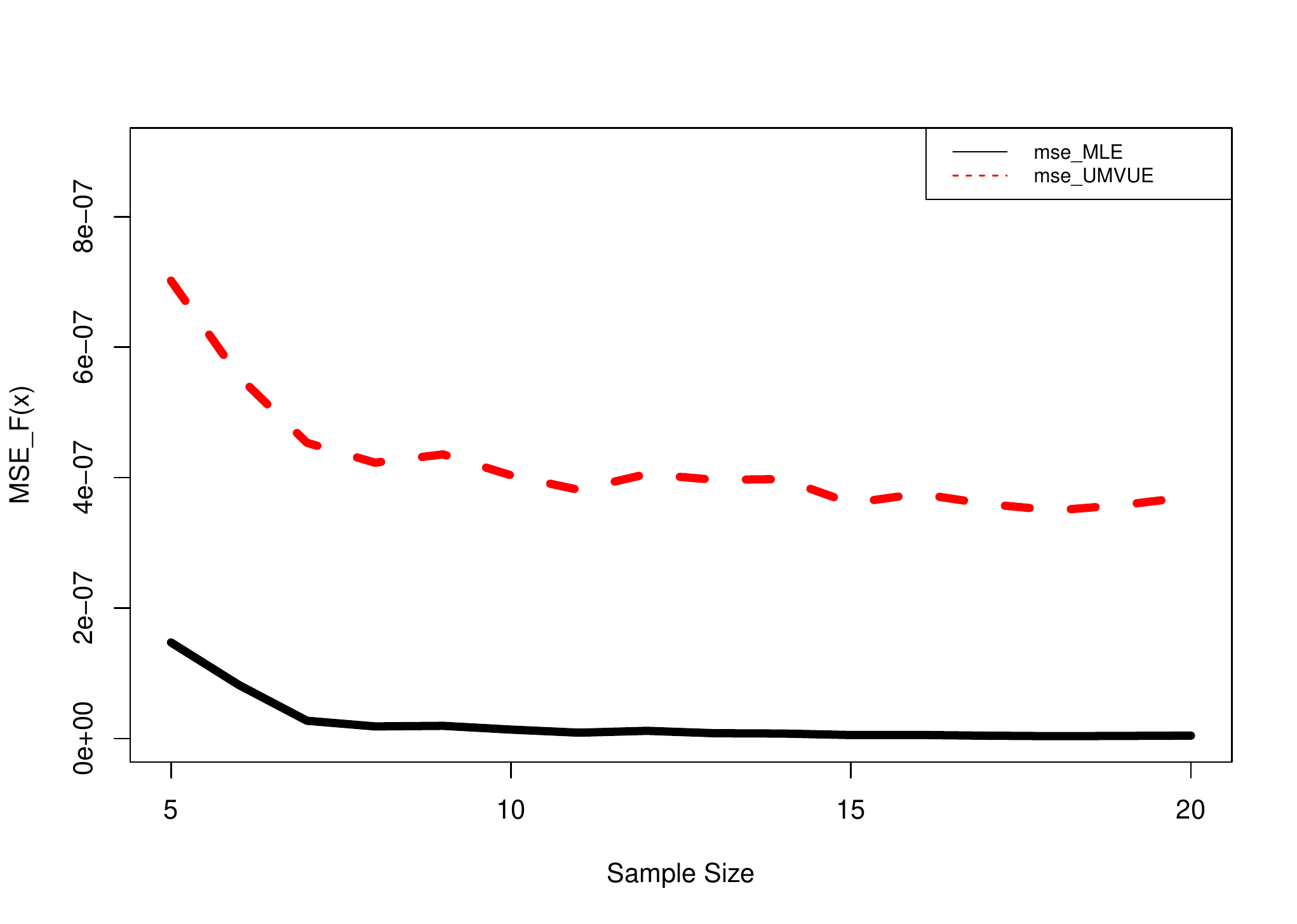}}
	\caption{ Graph of simulated MSE of the MLE and UMVUE of the PMF and the CDF of a NDL distribution for $x=2$, $\theta=0.01$ and $r=2$.}\label{ch6simpmfcdf}
\end{figure}

 A simulation study is carried out with $1,000(N)$ repetitions. Here we choose $a_0=1,~a_1=1$, $\theta = 0.01$, $x=2$ and $r = 2$. We compute MSE of the MLE and that of the UMVUE of the PMF and the CDF of a NDL distribution. From Figure $\ref{ch6simpmfcdf}$, it is clear that MSE decreases with an increasing sample size that shows the consistency property of the estimators. 
%In Figure $\ref{ch6simrelhazmrl}$, we have shown the simulated results of reliability, hazard rate and mean remaining life function using the PMF and the CDF of a NDL distribution.

\section{Data Analysis}\label{ch6data}
\begin{figure}[H]
	\centering
	\subfigure{\includegraphics[height=3in,width=2.6in]{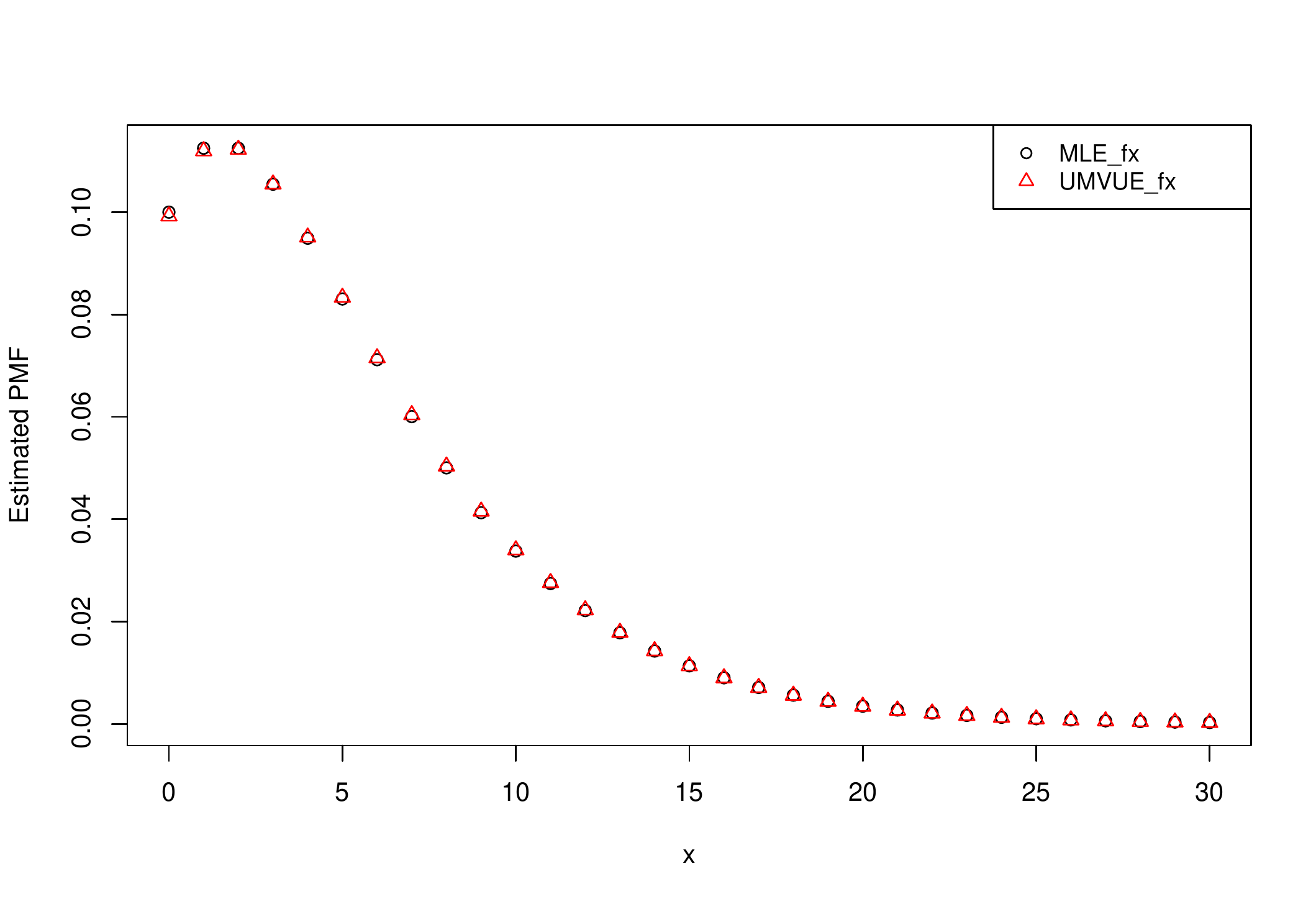}}
	\subfigure{\includegraphics[height=3in,width=2.6in]{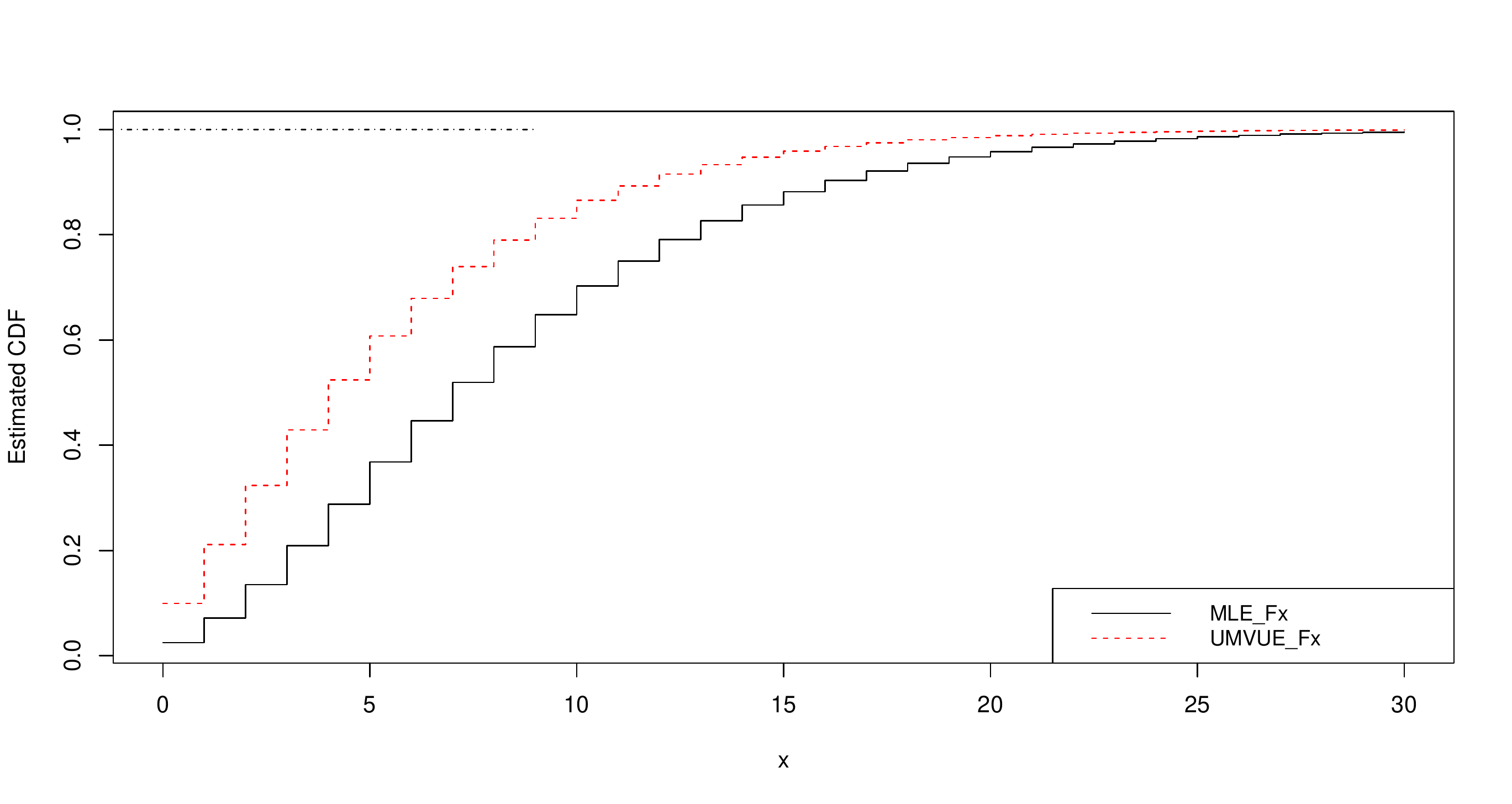}}
	\caption{Graph of the estimated PMF and that of the CDF of a NDL distribution.}\label{ch6datapmfcdf}
\end{figure}

\ We have studied the data comprise of numbers of fires in forest districts of Greece from period 1 July 1998 to 31 August 1998. The total number of observed samples is 123. This data set is obtained from Bakouch et al. (2014) and is shown in Table $\ref{ch6dataset}$. In Figure $\ref{ch6datapmfcdf}$, we have shown the estimated PMF and that of the CDF of a NDL distribution.
\begin{table}[H]
	\caption{\label{ch6dataset} Numbers of fires in Greece }
	\begin{center}
		{
			\footnotesize 
			\begin{tabular}{c c c c c c c c c c c c c c c c c c c}
				\hline
				2& 4& 4& 3& 3 &1& 2 &4& 3& 1& 1& 0& 5 &5& 0& 3& 1 &1 &0\\
				1& 0& 2& 0& 1& 2& 0& 0& 0& 0 &0 &0& 0& 0&0 &1& 4& 2& 2\\
				1 &2& 1& 2& 0 &2& 2& 1& 0 &3& 2& 1& 2& 2 &7& 3 &5& 2& 5\\
				4& 5& 6& 5 &4& 3 &8& 43 &8 &4 &4& 3& 10& 5& 4& 5& 12& 3 &8 \\
				12& 10 &11& 6& 1 &8& 9& 12& 9& 4& 8 &12& 11& 8& 6& 4& 7 &9& 15\\
				12& 15 &15 &12 &9 &16& 7& 11& 9& 11& 6 &5& 20& 9& 8 &8& 5& 7 &10\\
				6& 6 &5& 5 &15 &6 &8& 5& 6& & & & & & & & & & \\\hline
			\end{tabular}
		}
	\end{center}
\end{table}
\begin{table}[H]
	\caption{\label{ch6model}Model selection criterion}
	\begin{center}
		{
			\footnotesize 
			\begin{tabular}{|c|c|c|}
				\hline
				%\multicolumn{2}{|c|}{Negative log-likelihood}\\\hline
				%Estimators & NDL distribution\\\hline
				%\multirow{2}{*}{Estimators}&Negative log-likelihood\\
				%&NDL distribution\\\hline
				&\multicolumn{2}{|c|}{Negative log-likelihood}\\\cline{2-3}
				Estimators & \multicolumn{2}{|c|}{NDL distribution}\\\hline
				MLE & \multicolumn{2}{|c|}{340.0195}\\\hline
				UMVUE & \multicolumn{2}{|c|}{340.1765}  \\\hline
			\end{tabular}
		}
	\end{center}
\end{table}
\par We use the estimate of the negative log-likelihood values for the model selection criterion. Lower the value of negative log-likelihood indicates the better fit. From Table $\ref{ch6model}$, we found that MLE is better than UMVUE in a negative log-likelihood sense. 

\section{Concluding Remarks}\label{ch6conclusion}
\ In this article, two methods of estimation of the PMF and the CDF of a NDOPPE distribution have been considered. The MLE and the UMVUE have been found out. Simulation study is performed to compare the performances of the proposed methods of estimation. 
We have studied the performance of the estimators of the PMF and the CDF of a NDL distribution as a particular case of a NDOPPE distribution.
From simulation study results, it is found that UMVUE is better than MLE for the PMF and MLE is better than UMVUE for the CDF in MSE sense. From the model selection criterion, it is noticed that MLE is better than UMVUE in negative log-likelihood (i.e. AIC) sense. We analyze a data set to to illustrate fitting of the proposed distribution.

\end{document}